\input ./style/arxiv-general.cfg
\documentclass[seceqn,MSNbibl,number,citesort,dvips]{arxbj}
\makeatletter
   \@ifpackageloaded{graphicx}{}{\usepackage{graphicx}}
\makeatother
\usepackage{accents}
\usepackage{upgreek}

\aid{0}
\volume{21}
\issue{3}
\pubyear{2015}
\firstpage{1697}
\lastpage{1718}
\doi{10.3150/14-BEJ618} 

\makeatletter
\newcommand{\rrvert}{\vert}
\newcommand{\rrVert}{\Vert}
\newcommand{\llvert}{\vert}
\newcommand{\llVert}{\Vert}
\newtheorem{theorem}{Theorem}[section]
\newtheorem{corollary}[theorem]{Corollary}
\newtheorem{lemma}[theorem]{Lemma}
\newtheorem{proposition}[theorem]{Proposition}

\newproclaim{definition}[theorem]{Definition}

\newremark{remark}[theorem]{Remark}
\newremark{remarks}[theorem]{Remarks}
\newremark{example}[theorem]{Example}

\newcommand{\eqref}[1]{(\ref{#1})}
\newcommand{\RR}{\mathbb{R}}
\newcommand{\PP}{\mathbb{P}}
\newcommand{\CC}{\mathbb{C}}
\newcommand{\EE}{\mathbb{E}}
\newcommand{\cB}{\mathcal{B}}
\newcommand{\cF}{\mathcal{F}}
\newcommand{\cA}{\mathcal{A}}
\newcommand{\cL}{\mathcal{L}}
\newcommand{\dd}{\mathrm{d}}

\newcommand{\abs}[1]{\llvert  #1 \rrvert }
\newcommand{\norm}[1]{\llVert  #1 \rrVert }

\newcommand{\bbr}{\mathbb{R}}
\newcommand{\bbp}{\mathbb{P}}
\newcommand{\bbe}{\mathbb{E}}

\makeatother

\begin{document}
\begin{frontmatter}

\title{A criterion for invariant measures of It\^o processes based on
the symbol}
\runtitle{A criterion for invariant measures}

\begin{aug}
\author[1]{\inits{A.}\fnms{Anita}~\snm{Behme}\corref{}\thanksref{1}\ead[label=e1]{behme@ma.tum.de}} \and
\author[2]{\inits{A.}\fnms{Alexander}~\snm{Schnurr}\thanksref{2}\ead[label=e2]{alexander.schnurr@math.tu-dortmund.de}}
\address[1]{Institut f\"ur Mathematische Stochastik, Technische
Universit\"at Dresden,
D-01062 Dresden, Germany
and Technische Universit\"at
M\"unchen, Center for Mathematical Sciences, Boltzmannstrasse 3,
D-85748 Garching, Germany.
\printead{e1}}
\address[2]{Department Mathematik, Universit\"at Siegen,
Walter-Flex-Str. 3, D-57072 Siegen and Technische Universit\"at
Dortmund, Lehrstuhl IV, Fakult\"at f\"ur Mathematik, D-44227 Dortmund,
Germany.\\\printead{e2}}
\end{aug}

\received{\smonth{10} \syear{2013}}

%
\begin{abstract}
An integral criterion for the existence of an invariant measure of an
It\^o process is developed. This new criterion is based on the
probabilistic symbol of the It\^o process. In contrast to the standard
integral criterion for invariant measures of Markov processes based on
the generator, no test functions and hence no information on the domain
of the generator is needed.
\end{abstract}

%
\end{frontmatter}

\section{Introduction}\label{S1}

Consider $(\Omega, \cF, (\cF_t)_{t\geq0}, (X_t)_{t\geq0}, \bbp
^x)_{x\in\bbr^d}$ to be a Feller process on $\RR^d$ with semigroup
$(T_t)_{t\geq0}$ on $C_0(\RR^d)$, that is,
\[
T_tf(x)=\int_{\RR^d} f(y)\mu_t(x,\dd
y)=\EE^x\bigl[f(X_t)\bigr],
\]
where $\mu_t(x,\dd y)=\PP^x(X_t\in\dd y)=\PP(X_t\in\dd y|X_0=x)$
are the transition probabilities and $C_0(\RR^d)$ are the continuous,
real-valued functions on $\RR^d$ vanishing at infinity.
Then the {\sl infinitesimal generator} $\cA$ of $X$ is defined by
\[
\cA f=\lim_{t\to0} \frac{T_tf-f}{t}
\]
for all functions $f$ in the domain of $\cA$, that is, all $f$ in
\[
D(\cA)= \biggl\{f\in C_0\bigl(\RR^d\bigr), \lim
_{t\to0} \frac{T_tf-f}{t} \mbox{ exists in } \|\cdot
\|_\infty \biggr\}.
\]
It is known (see, e.g., \cite{liggett}, Theorem~3.37) that a probability
measure $\mu$ on $\RR^d$ is {\sl invariant} (or {\sl stationary}) for
the Feller process $X$ with semigroup $(T_t)_{t\geq0}$, meaning that
\[
\int_{\RR^d} T_tf(x) \,\dd\mu(x)= \int
_{\RR^d} f(x)\,\dd\mu (x),\quad\quad \forall f\in C_0\bigl(
\RR^d\bigr), t\geq0,
\]
if and only if
\begin{equation}
\label{generatorstat} \int_{\RR^d} \cA f(x)\mu(\dd x) = 0
\end{equation}
holds for all $f$ in a core of the generator $\cA$.

In the special case that $X$ is a rich Feller process, that is, a Feller
process with the property that the test functions $C_c^\infty(\RR^d)$
are contained in the domain $D(\cA)$ of the generator, this generator
is a pseudo-differential operator with negative definite symbol
$p(x,\xi)$
\begin{equation}
\label{generatorsymbol} \cA f(x)=-\int_{\RR^d} \mathrm{e}^{\mathrm{i}x'\xi}
p(x,\xi) \hat{f}(\xi) \,\dd\xi
\end{equation}
for all $f\in C_c^\infty(\RR^d)$ (see, e.g., \cite{ReneLM13}, Definition~2.25 and
Corollary~2.23). Hereby the superscript ``$\,'$'' denotes the transpose
of a vector, $\hat{f}(y)=(2\uppi)^{-d}\int_{\RR^d} \mathrm{e}^{-\mathrm{i}x'y} f(x) \,\dd
x$ denotes the Fourier transform of $f$ and $C_c^\infty(\RR^d)$ is
the space of infinitely often continuously differentiable functions on
$\RR^d$ with compact support.
Thus in this case, equations
\eqref{generatorstat} and \eqref{generatorsymbol} together yield that
if $\mu$
is an invariant law for $X$, then for all $f$ in $C_c^\infty(\RR^d)$
one has
\begin{equation}
\label{symbolstatfeller} \int_{\RR^d} \int_{\RR^d}
\mathrm{e}^{\mathrm{i}x'\xi} p(x,\xi) \hat {f}(\xi) \,\dd\xi\mu(\dd x) =0.
\end{equation}
Conversely, if \eqref{symbolstatfeller} holds
for all $f$ in a core $D\subset C_c^\infty(\RR^d)\cap D(\cA)$ of
$\cA$, then $\mu$
is an invariant law for $X$.

Further, formally applying Fubini's theorem on equation \eqref
{symbolstatfeller} leads to
\begin{equation}
\label{symbolstat} \int_{\RR^d} \mathrm{e}^{\mathrm{i}x \xi} p(x,\xi)
\mu(\dd x)=0,\quad\quad \mbox{for } \lambda\mbox{-almost all } \xi\in\RR^d,
\end{equation}
where $\lambda$ denotes the Lebesgue measure. Observe that \eqref
{symbolstat} is an equation directly relating the symbol to the
invariant law and does not involve any test functions. This is a big
advantage for application of \eqref{symbolstat} compared to \eqref
{generatorstat} where we started from.

As a trivial example of application, consider $L$ to be a L\'evy
process, that is, we have $p(x,\xi)=\psi_L(\xi)$ where $\psi_L$
denotes the L\'evy--Khintchine exponent of $L$. Then \eqref{symbolstat}
is equivalent to $\psi_L(\xi)\phi_\mu(\xi)=0$ where $\phi_\mu$
denotes the characteristic function of $\mu$. Since this
characteristic function is continuous and takes the value $1$ at $0$
this yields that $\psi_L(\xi)=0$ for $\xi$ in some neighborhood of
$0$ which implies that $\psi_L(\xi)=0$ for all $\xi$ (cf. \cite
{sato}, Lemma~13.9). Hence, $L_t=0$. Indeed the zero process is the
only L\'evy process which admits an invariant law (cf. \cite{sato}, Exercise~19.6). On the contrary, for the zero process \eqref{symbolstat}
follows immediately.

Instead of making the above computations rigorous in the case of rich
Feller processes, in this paper we will consider a wider class of
processes. Therefore recall that in the case of general Markov
processes, necessity of equation \eqref{generatorstat} for $\mu$ to
be invariant is still given. For example, \cite{ethierkurtz}, Proposition~9.2, states that
if a distribution $\mu$ is invariant for a Markov process $X$ then
\eqref{generatorstat} holds for all $f$ in the domain of the generator
of $X$.
Remark that one part of the literature on Markov processes (and so
\cite{ethierkurtz}) defines the generator on functions with bounded
support, that is, in $C_b(\RR^d) \supset C_0(\RR^d)$ which
does not fit into the setting we described above and which we will use
throughout this paper.

For the converse direction in the Markovian setting, that is, to show
sufficiency of \eqref{generatorstat} for $\mu$ to be invariant,
further assumptions on the generator are needed as discussed for
example, in \cite{echeverria} and \cite{bhatt}.
This is the reason why sometimes (e.g., in \cite{alberuediwu}) a
probability measure $\mu$ is called {\sl infinitesimal invariant} for
a given generator and domain, if and only if \eqref{generatorstat} is
fulfilled for all $f$ in the domain of the generator.

General Markov processes do not necessarily have an associated symbol.
Hence in this paper, we restrict ourselves to It\^o processes as they
are defined below (Definition~\ref{def:ito}). This class includes the rich
Feller processes, but is much more general. For these It\^o processes,
we derive the relation between the symbol of an It\^o process as
defined below and its invariant law. In particular, our aim is to show
in as much generality as possible, that for an It\^o process $X$ with
symbol $p(x,\xi)$ equation \eqref{symbolstat} holds if (and only if)
$\mu$ is an (infinitesimal) invariant law for $X$.

The paper is outlined as follows. In Section~\ref{S2}, we recall the
necessary definitions of It\^o processes and symbols as they will be
used throughout this paper. Section~\ref{S3} then shows necessity of
\eqref{symbolstat} for $\mu$ to be an invariant law for a wide class
of It\^o processes. Several examples are given and some special cases
are studied. Sufficiency of \eqref{symbolstat} for $\mu$ to be
infinitesimal invariant is then treated in Section~\ref{S4} and again
it is illustrated by special cases. Some rather technical proofs for
results in Section~\ref{S3} have been postponed to the closing Section~\ref{Sproof}.

\section{Preliminaries}\label{S2}

In 1998, Jacob came up with the idea to use a probabilistic approach in
order to calculate the so-called symbol of a stochastic process \cite
{Jaco1998}. This probabilistic formula was generalized in the same year
to rich Feller processes by Schilling \cite{schilling98pos}. Let us
recall the definition.

\begin{definition}
Let $(X_t)_{t\geq0}$ be a Markov process in $\RR^d$. Define for every
$x,\xi\in\bbr^d$ and $t\geq0$ the quantity
\[
\lambda_\xi(x,t):= - \frac{\bbe^x \mathrm{e}^{\mathrm{i}(X_t-x)'\xi}-1}{t}.
\]
We call $p\dvtx \bbr^d\times\bbr^d\to\CC$ given by
\begin{equation}
\label{symbollimit} p(x,\xi):=-\lim_{t\downarrow0}\frac{\bbe^x \mathrm
{e}^{\mathrm{i}(X_t-x)'\xi}-1}{t}=\lim
_{t\downarrow0} \lambda_\xi(x,t)
\end{equation}
the \emph{probabilistic symbol} of $X$ if the limit exists for every
$x,\xi$.
\end{definition}

For rich Feller processes satisfying the growth condition
\begin{equation}
\label{growth} \sup_{x\in\bbr^d} \bigl|p(x,\xi)\bigr|\leq c\bigl(1+\llVert \xi
\rrVert ^2\bigr), \quad\quad\xi\in\bbr^d,
\end{equation}
for $\|\cdot\|$ denoting an arbitrary submultiplicative norm, the
probabilistic symbol and the symbol in equation \eqref
{generatorsymbol} coincide. This is why we use the letter $p$ and call
the new object again ``symbol''. Remark that \eqref{growth} is a
standard condition in this context. We refer to \cite{ReneLM13} for a
comprehensive overview on Feller processes and their symbol.

The class of rich Feller processes includes L\'evy processes as special
case. For these, the symbol only depends on $\xi$ and coincides with
the {\em L\'evy exponent}, that is,
\[
p(x,\xi)=\psi_L(\xi):=-\log\EE\bigl[\mathrm{e}^{\mathrm{i}L_1'\xi}\bigr]
=-\mathrm{i}
\ell' \xi+ \frac{1}2 \xi' Q \xi- \int
_{\RR^d} \bigl(\mathrm{e}^{\mathrm{i}\xi'y} -1 - \mathrm{i}\xi
' y 1_{\{\|y\|<1\}}(y) \bigr) N(\dd y),
\]
where $(L_t)_{t\geq0}$ is a L\'evy process with {\em characteristic
triplet} $(\ell, Q, N)$. For details on L\'evy processes in
particular, we refer to \cite{sato}.

On the other hand, every rich Feller process is an It\^o process (cf.
\cite{symbolkillingopenset}, Theorem~3.9) in the following sense. It is
this class we are dealing with in the present paper.

\begin{definition} \label{def:ito}
An \emph{It\^o process} is a strong Markov process, which is a
semimartingale w.r.t. every $\bbp^x$ having semimartingale
characteristics of the form
\begin{eqnarray}
B_t^{(j)}(\omega) &=&\int_0^t
\ell^{(j)}\bigl(X_s(\omega)\bigr) \,\dd s,\quad\quad j=1,\ldots,d\nonumber
\\
C_t^{jk}(\omega) &=&\int_0^t
Q^{jk}\bigl(X_s(\omega)\bigr) \,\dd s,\quad\quad j,k=1,\ldots,d
\\
\nu(\omega;\dd s,\dd y) &=&N\bigl(X_s(\omega),\dd y\bigr) \,\dd s\nonumber
\end{eqnarray}
for every $x\in\bbr^d$ with respect to a fixed cut-off function $\chi
$. Here, $\ell(x)=(\ell^{(1)}(x),\ldots,\ell^{(d)}(x))'$ is a vector in
$\bbr^d$, $Q(x)$ is a positive semi-definite matrix and $N$ is a Borel
transition kernel such that $N(x,\{0\})=0$. We call $\ell$, $Q$ and
$n:=\int_{y\neq0} (1\wedge\norm{y}^2) N(\cdot,\dd y)$ the \emph
{differential characteristics} of the process.
\end{definition}

Usually we will have to impose the following condition on the
differential characteristics of It\^o processes.

\begin{definition} \label{def:finelycont}
Let $X$ be a Markov process and $f\dvtx \bbr^d\to\bbr$ be a
Borel-measurable function. Then $f$ is called $X$\emph{-finely
continuous} (or \emph{finely continuous}, for short) if the function
\begin{equation}
\label{rightcont} t\mapsto f(X_t)=f\circ X_t
\end{equation}
is right continuous at zero $\bbp^x$-a.s. for every $x\in\bbr^d$.
\end{definition}

\begin{remark}
Fine continuity is usually introduced in a different way
(see \cite{blumenthalget}, Section II.4, and~\cite{fuglede}). However, by \cite
{blumenthalget}, Theorem~4.8, this is equivalent to \eqref{rightcont}. It
is this kind of right continuity which we will use in our proofs. Let
us mention that this assumption is very weak, even weaker than ordinary
continuity.
\end{remark}

In \cite{cinlarjacod81}, the class of It\^o processes in the sense of
Definition~\ref{def:ito} has been characterized as the set of
solutions of very general SDEs. In particular, as mentioned already,
the class of It\^o processes contains the class of rich Feller
processes. The following example (cf. \cite{detmp1}, Example~5.2)
shows that this inclusion is strict: The process
\[
X_t^x=\cases{ x-t & \quad$\mbox{under } \bbp^x
\mbox{ for } x<0$,
\cr
0 &\quad$\mbox{under } \bbp^x \mbox{ for } x=0$,
\cr
x+t &\quad$\mbox{under } \bbp^x \mbox{ for } x>0$,}\quad\quad  t\geq0,
\]
is an It\^o process with bounded and finely continuous differential
characteristics which is not Feller.

For It\^o processes, we can compute the probabilistic symbol. In
particular, we even have the following connection of probabilistic
symbol and generator.

\begin{lemma} \label{itogenerator}
If the test functions $C_c^\infty(\bbr^d)$ are contained in the
domain $D(\cA)$ of the generator $\cA$ of an It\^o process $X$, the
representation
\eqref{generatorsymbol} holds for every $f\in C_c^\infty(\bbr^d)$,
where $p(x,\xi)$ is the probabilistic symbol.
\end{lemma}

\begin{pf}
An operator $\overline{\cA}$ with domain $D(\overline{\cA})$ is
called the extended generator of $X$, if $D(\overline{\cA})$ consists
of those Borel measurable functions $f$ for which there exists a $(\cB
^d)^*$-measurable function $\overline{\cA}f$ such that the process
\[
M_t^f:=f(X_t)-f(X_0)-\int
_0^t \overline{\cA}f(X_s) \,\dd s
\]
is a local martingale (cf. \cite{vierleute}, Definition~7.1). Here, $(\cB
^d)^*$ denotes as usual the universally measurable sets (see, e.g.,
\cite{blumenthalget}, Section~0.1). We have $C_c^\infty(\bbr
^d)\subset D(\cA)\subset D(\overline{\cA})$ and $\overline{\cA
}|(D(\cA))=\cA$ by Dynkin's formula. By \cite{vierleute}, Theorem~7.16,
we obtain that
\begin{eqnarray*}
-\overline{\cA}f(x)&=&-\sum_{j\leq d}
\ell^{(j)}(x) D_jf(x)-\frac
{1}{2} \sum
_{j,k\leq d} Q^{j,k}D_{j,k}f(x)
\\
&&{}-\int_{y\neq0} f(x+y)-f(x)-\chi(y)\sum
_{j\leq d} y^{(j)}D_jf(x) N(x,\dd y)
\\
&=&\int_{\bbr^d} \mathrm{e}^{\mathrm{i}x'\xi} p(x,\xi) \widehat{f}(
\xi) \,\dd\xi
\end{eqnarray*}
for every $f\in C_c^\infty(\bbr^d)$ which is \eqref{generatorsymbol}.
\end{pf}

To end this section, let us also mention that there exists a formula to
calculate the symbol even for the wider class of homogeneous diffusions
with jumps in the sense of \cite{jacodshir} (cf. \cite{alexgeneralsymbol}, Theorem~3.6). However, this formula uses stopping times and
can not be used in our considerations. In the proof of Theorem~\ref
{thm:nec} below, we have to use the classical version of the
probabilistic symbol presented above.

For the even wider class of Hunt semimartingales, the limit \eqref
{symbollimit} is not defined and hence the symbol does not exist any
more (cf. \cite{manstaviciusschnurr}).

\section{Necessity}\label{S3}

We start by showing the necessity of \eqref{symbolstat} for $\mu$ to
be an invariant law.

\begin{theorem} \label{thmnectest} Let $(X_t)_{t\geq0}$ be an It\^o
process with generator $\cA$ whose domain $D(\cA)$ contains the test
functions $C_c^\infty(\bbr^d)$ and with symbol $p(x,\xi)$. Assume
$\mu$ is an invariant measure for $X$ such that $\int_{\RR^d}
|p(x,\xi)| \mu(\dd x)<\infty$.
Then
\[
\int_{\RR^d} \mathrm{e}^{\mathrm{i}x'\xi} p(x,\xi) \mu(\dd x)=0
\quad\quad\forall\xi\in\RR^d.
\]
\end{theorem}

\begin{pf}
It is well known, that for an invariant measure $\mu$ it holds
\[
\int_{\RR^d} \cA f(x)\mu(\dd x)=0 \quad\quad\mbox{for all } f \in D(\cA).
\]
By Lemma~\ref{itogenerator}, the generator $\cA$ admits the
representation \eqref{generatorsymbol} for all $f\in C_c^\infty(\RR
^d)$. Hence using Fubini's theorem, we obtain for all $f\in C_c^\infty
(\RR^d)$ 
\begin{eqnarray*}
0 &=&\int_{\RR^d} \cA f(x)\mu(\dd x)= - \int
_{\RR^d} \int_{\RR^d} \mathrm{e}^{\mathrm{i}x'\xi}
p(x,\xi) \hat{f}(\xi) \,\dd\xi\mu(\dd x)
\\
&=& - \int_{\RR^d} \hat{f}(\xi) \int_{\RR^d}
\mathrm{e}^{\mathrm{i}x'\xi} p(x,\xi) \mu(\dd x) \,\dd\xi.
\end{eqnarray*}
This yields that for $\lambda$-a.a. $\xi$ it holds $\int_{\RR^d}
\mathrm{e}^{\mathrm{i}x'\xi} p(x,\xi) \mu(\dd x)=0$. Since by hypothesis
$p(x,\xi)$ is absolutely integrable with respect to $\mu$, its
Fourier transform with respect to $\mu$ is continuous. This gives the claim.
\end{pf}

\begin{example} \label{exampleGOUnec}
Let $(X_t)_{t\geq0}$ be a generalized Ornstein--Uhlenbeck process,
defined as the unique solution of
\[
\dd X_t=X_{t-} \,\dd U_t + \dd
L_t,\quad\quad t\geq0,
\]
for two independent L\'evy processes $(U_t)_{t\geq0}$ and
$(L_t)_{t\geq0}$. It has been shown in \cite{behmelindner13}, Theorem~3.1, that $X$ is a Feller process, that the domain of
its generator contains $C_c^\infty(\bbr)$ and that $C_c^\infty(\bbr
)$ is a core for the generator.

Further it follows from the results
in \cite{sdesymbol} that the symbol of $X$ is given by
\[
p(x,\xi)= \psi_U(x\xi)+\psi_L(\xi), \quad\quad x, \xi\in\RR.
\]
Assume $\mu$ is a probability measure on $\RR$ such that $\int x^2
\mu(\dd x)<\infty$. Then due to the specific form of the symbol $\int_{\RR^d} |p(x,\xi)| \mu(\dd x)<\infty$ is automatically fulfilled.
Hence, we see from Theorem~\ref{thmnectest} that
\begin{equation}
\label{eqGOUstat} \int_{\RR} \mathrm{e}^{\mathrm{i}x\xi}
\psi_U(x\xi) \mu(\dd x)=
-\psi_L(\xi)
\phi_\mu(\xi)
\end{equation}
is a necessary condition for $\mu$ to be invariant for $X$. Remark
that equation \eqref{eqGOUstat} has also been obtained in
\cite{behmelindner13}, Theorem~4.1.
\end{example}

In general, we have only little or no information on the domain of the
generator of an It\^o process which makes the above theorem
inapplicable. We will see in the following, that it is possible to
obtain similar results without any information on the domain of the
generator by using the probabilistic definition of the symbol directly.
The first case we consider is the It\^o process with bounded and finely
continuous differential characteristics. Examples include Feller
processes satisfying the growth condition \eqref{growth}. Although the
boundedness assumption seems to be rather restrictive, this class of
processes already contains various interesting examples, which are used
in stochastic modeling and mathematical statistics (see, e.g., Example~\ref{ex:bounded} below).

\begin{theorem} \label{thm:nec}Let $(X_t)_{t\geq0}$ be an $\RR
^d$-valued It\^o process with bounded, finely continuous differential
characteristics which admits an invariant law $\mu$ and
whose symbol is given by $p(x,\xi)$, $x,\xi\in\RR^d$. Then
\[
\int_{\RR^d} \mathrm{e}^{\mathrm{i}x' \xi} p(x,\xi) \mu(
\dd x)=0 \quad\quad\mbox{for all } \xi\in\RR^d.
\]
\end{theorem}

For the proof of Theorem~\ref{thm:nec}, we need the following lemma,
which also shows the form of the symbol in the given setting.

\begin{lemma} \label{glmkonv}
Let $(X_t)_{t\geq0}$ be an $\RR^d$-valued It\^o process with bounded,
finely continuous differential characteristics.
For every $\xi\in\bbr^d$ the limit
\[
p(x,\xi):= -\lim_{t\downarrow0}\frac{\bbe^x \mathrm
{e}^{\mathrm{i}(X_t-x)'\xi}-1}{t}=\lim
_{t\downarrow0} \lambda_\xi(x,t)
\]
exists and the functions $\lambda_\xi$ are globally bounded in $x$
(and $t$) for every $\xi\in\bbr^d$. As the limit, we obtain
\begin{equation}
\label{symbol} p(x,\xi)=-\mathrm{i}\ell(x)'\xi+ \frac{1}{2}
\xi'Q(x) \xi-\int_{y\neq0} \bigl(
\mathrm{e}^{\mathrm{i}y'\xi}-1 -\mathrm{i}y'\xi\cdot\chi(y) \bigr) N(x,\dd y).
\end{equation}
\end{lemma}

The proof of Lemma~\ref{glmkonv} is postponed to Section~\ref{Sproof}.

\begin{pf*}{Proof of Theorem~\ref{thm:nec}}
Let $X_\infty$ be a random variable such that $\mu=\cL(X_\infty)$
where ``$\cL$'' stands for ``law of''. Then using Lemma~\ref
{glmkonv}, we obtain by Lebesgue's dominated convergence theorem
\begin{eqnarray*}
\int_{\RR^d} \mathrm{e}^{\mathrm{i}x' \xi} p(x,\xi) \mu(\dd x) &=&
-\int_{\RR^d} \mathrm{e}^{\mathrm{i}x' \xi} \lim
_{t\to0} \EE^x \biggl[\frac{\mathrm{e}^{\mathrm{i}(X_t-x)'\xi}-1}{t} \biggr] \mu(
\dd x)
\\
&=& -\lim_{t\to0}\frac{1}{t} \int_{\RR^d}
\mathrm{e}^{\mathrm{i}x'\xi} \EE^x \bigl[\mathrm{e}^{\mathrm{i}(X_t-x)'\xi}-1
\bigr] \mu(\dd x),
\end{eqnarray*}
with
\begin{eqnarray*}
\int_{\RR^d} \mathrm{e}^{\mathrm{i}x' \xi} \EE^x \bigl[
\mathrm{e}^{\mathrm{i}(X_t
-x)' \xi}-1 \bigr] \mu(\dd x) &=& \int_{\RR^d}
\EE^x \bigl[\mathrm{e}^{\mathrm{i}X_t' \xi} \bigr] \mu (\dd x) - \int
_{\RR^d} \mathrm{e}^{\mathrm{i}x' \xi} \mu(\dd x)
\\
&=& \int_{\RR^d} \EE \bigl[\mathrm{e}^{\mathrm{i}X_t'\xi}|X_0=x
\bigr] \mu(\dd x) - \EE \bigl[ \mathrm{e}^{\mathrm{i} X_\infty' \xi} \bigr]
\\
&=& \EE \bigl[ \mathrm{e}^{\mathrm{i} X_\infty' \xi} \bigr]- \EE \bigl[
\mathrm{e}^{\mathrm{i} X_\infty' \xi} \bigr]
\\
&=&0.
\end{eqnarray*}
\upqed\end{pf*}

The following example is taken from \cite{lamanga}, Section~5.7. It
is derived by a transformation from a classical example due to
Barndorff-Nielsen \cite{bn78}.

\begin{example} 
\label{ex:bounded}
Let $(X_t)_{t\geq0}$ be the unique solution of the SDE
\[
\dd X_t= b(X_t) \,\dd t + \sigma(X_t) \,\dd
W_t,\quad\quad t\geq0,
\]
with $X_0=x_0$, a standard Brownian motion $(W_t)_{t\geq0}$,
\[
b(x)= - \biggl(\vartheta+\frac{c^2}{2 \cosh(x)} \biggr) \frac{\sinh
(x)}{\cosh^2(x)} \quad\mbox{and}\quad \sigma(x)=\frac{c}{\cosh(x)},
\]
where $\vartheta,c>0$. For $x_0\in\bbr$ the \emph{scale density}
and the \emph{speed density} of $X$ are then given by
\[
s(x):=\exp \biggl(-2 \int_{x_0}^x
\frac{b(u)}{\sigma^2(x)} \,\dd u \biggr) \quad\mbox{and}\quad
 m(x):=\frac{1}{\sigma^2(x)s(x)}.
\]
Since $\int s(x) \,\dd x =\infty$ while $M:= \int m(x) \,\dd x <
\infty$, we are in the setting of \cite{lamanga}, Section~5.2.
There, the authors restate the well-known fact that the unique
stationary distribution of the process $X$ in this case admits the density
\[
\pi(x)=\frac{m(x)}{M}.
\]
By our above result, this means
\[
\int\mathrm{e}^{\mathrm{i}x\xi} \bigl(\bigl\llvert \sigma(x)\bigr\rrvert
^2 \llvert \xi\rrvert ^2 - \mathrm{i}b(x)\xi \bigr) \pi(x) \,\dd x
=0
\]
since $p(x,\xi)=\abs{\sigma(x)}^2 \abs{\xi}^2 - \mathrm{i}b(x)\xi$ is the
symbol of $X$ by \cite{sdesymbol}, Theorem~3.1.
\end{example}

\subsection{L\'evy driven SDEs}

In general, we can not drop the boundedness assumption on the
differential characteristics which we have used in Theorem~\ref
{thm:nec}. This assumption corresponds to bounded coefficients of the
SDEs whose solutions are the considered It\^o processes. However, in
some cases we are able to generalize our result as shown in the
following proposition where a linearly growing coefficient is allowed.
Another possible extension of Theorem~\ref{thm:nec} will be stated in
Proposition~\ref{corbrown} below.

\begin{proposition} \label{proplevy}
Let $(X_t)_{t\geq0}$ be the unique solution of the SDE
\[
\dd X_t=-aX_t \,\dd t + \Phi(X_{t-}) \,\dd
L_t,\quad\quad t\geq0,
\]
where $a\in\RR$, $\Phi\dvtx \bbr^d\to\bbr^{d\times n}$ is bounded and
locally Lipschitz continuous and $(L_t)_{t\geq0}$ is a L\'evy process
in $\RR^n$ satisfying $\EE\|L_1\|<\infty$.
Then $X$ is an It\^o process and for every $\xi\in\bbr^d$ the limit
$p(x,\xi)=\lim_{t\downarrow0} \lambda_\xi(x,t)$
exists and the functions $\lambda_\xi$ are globally bounded in $x$
(and $t$) for every $\xi\in\bbr^d$. Furthermore if $\mu$ is an
invariant law of $X$, then
\[
\int_{\RR^d} \mathrm{e}^{\mathrm{i}x' \xi} p(x,\xi) \mu(
\dd x)= \int_{\RR
^d} \mathrm{e}^{\mathrm{i}x' \xi}\bigl(
\psi_L\bigl(\Phi(x)'\xi\bigr)+\mathrm{i}ax'\xi
\bigr) \mu(\dd x) =0\quad\quad \mbox{for all } \xi\in\RR^d.
\]
\end{proposition}

The proof of Proposition~\ref{proplevy} is postponed to Section~\ref{Sproof}.

\begin{remarks}
\begin{enumerate}[(iii)]
\item[(i)] The structure of the symbol in Proposition~\ref{proplevy} is not
surprising. It is exactly what one would get for the generalized symbol
which uses stopping times (see \cite{alexgeneralsymbol}, Theorem~3.6).
However, it is important to see that for the classic probabilistic
symbol, as we have to use it in our context, the convergence of the
$\lambda_\xi(x,t)$ is uniform.
\item[(ii)] The class of It\^o processes studied in the above proposition is
a subset of the class of solutions of delay equations for which
stationarity was treated in \cite{reiss}. Remark that even in their
general paper, the authors have to impose the boundedness condition on
$\Phi$ (\cite{reiss}, Assumption~4.1(c)).
\item[(iii)] Although the assumption of $L$ having a finite first moment seems
very restrictive, it cannot be released.

As an example consider the L\'evy-driven Ornstein--Uhlenbeck (OU) process
\[
X_t=\mathrm{e}^{-\lambda t} \biggl(X_0+\int
_{(0,t]}\mathrm {e}^{\lambda s}\,\dd L_s
\biggr),\quad\quad t\geq0,
\]
for $\lambda>0$ and a symmetric $\alpha$-stable, real-valued L\'evy
process $(L_t)_{t\geq0}$, $\alpha\in[1,2]$. Then $(X_t)_{t\geq0}$
solves the SDE $\dd X_t=-\lambda X_{t-} \,\dd t + \dd L_t$ and by results
in \cite{sdesymbol} its symbol is given by
\[
p(x,\xi)= \mathrm{i} \lambda x\xi+ \psi_L(\xi) = \mathrm{i} \lambda x\xi+ |\xi
|^\alpha.
\]
Since $L$ is $\alpha$-stable, we have $E|L_1|^r<\infty$ if and only
if $r<\alpha$. Hence in particular $E\log^+ |L_1|<\infty$, such that
by \cite{lindnermaller}, Theorem~2.1, the process
$(X_t)_{t\geq0}$ admits a stationary solution with distribution $\cL
(X_\infty)$. By \cite{behme2011}, Theorem~3.1, it further follows that
$E|X_\infty|^r<\infty$ for all $r<\alpha$ and~-- as it was already
remarked in \cite{behme2011}~-- this result is sharp. Hence $\cL
(X_\infty)$ does not necessarily have a finite $\alpha$th moment and
the integral \eqref{symbolstat} does not necessarily exist.
\end{enumerate}
\end{remarks}

In \cite{alberuediwu}, the authors studied the absolutely continuous
invariant measures of solutions of SDEs of the form
\begin{equation}
\label{alberuediwuSDE} \dd X_t=\sqrt{2a_1} \,\dd
W_t+ \beta(X_t) \,\dd t + a_2 \,\dd
Z_t,\quad\quad t\geq0,
\end{equation}
for suitable coefficients, a standard Brownian motion $(W_t)_{t\geq0}$
and a pure-jump process $(Z_t)_{t\geq0}$ with stable (type) L\'evy measure.

For such SDEs with additive L\'evy noise, we obtain the following
proposition whose proof can be found in Section~\ref{Sproof}. Remark
that this proposition would follow directly from Proposition~\ref
{proplevy}, if we imposed $\EE\|L_1\|<\infty$ and $\EE\|Z_1\|<\infty$.

\begin{proposition} \label{proplevytwo} Let $(X_t)_{t\geq0}$ be the
unique solution of the SDE
\[
\dd X_t= b \,\dd Z_t + \Phi(X_{t-}) \,\dd
L_t, \quad\quad t\geq0,
\]
where $(L_t)_{t\geq0}$ and $(Z_t)_{t\geq0}$ are independent L\'evy
processes which are $n$- respectively $d$-dimensional, $b\in\bbr$ and
$\Phi\dvtx  \bbr^d\to\bbr^{d\times n}$ is bounded and locally Lipschitz.
Then $X$ is an It\^o process, for every $\xi\in\bbr^d$ the limit
$p(x,\xi)=\lim_{t\downarrow0} \lambda_\xi(x,t)$
exists and the functions $\lambda_\xi$ are globally bounded in $x$
(and $t$) for every $\xi\in\bbr^d$. Furthermore if $\mu$ is an
invariant law of $X$, then
\[
\int_{\RR^d} \mathrm{e}^{\mathrm{i}x' \xi} p(x,\xi) \mu(
\dd x)=\int_{\RR
^d} \mathrm{e}^{\mathrm{i}x' \xi} \bigl(
\psi_L\bigl(\Phi(x)'\xi\bigr)+ \psi_Z(b\xi
) \bigr) \mu(\dd x) =0 \quad\quad\mbox{for all } \xi\in\RR^d.
\]
\end{proposition}

Further, restricting to the class of processes whose absolutely
continuous invariant measures have been studied in \cite{alberuediwu},
we obtain the following corollary.

\begin{corollary} \label{cor:albeverio}
Let $(X_t)_{t\geq0}$ be the unique solution of the SDE \eqref{alberuediwuSDE}
where $a_1, a_2 \geq0$, $a_1+a_2>0$, $\beta\dvtx \RR^d\to\RR^d$ is
Borel measurable, locally Lipschitz, bounded and its Fourier transform
exists, $(W_t)_{t\geq0}$ is an $\RR^d$-valued standard Brownian
motion, and $(Z_t)_{t\geq0}$ is an $\RR^d$-valued pure-jump process
with L\'evy measure $\nu_\alpha(\dd y):=|y|^{-(d+\alpha)} \,\dd y$,
for $\alpha\in(0,2)$.

Assume $\mu(\dd x)=\rho(x)\,\dd x$ is invariant for $X$. Then
\begin{equation}
\label{albecriterion} \bigl(a_1|\xi|^2- a_2
c_\alpha|\xi|^\alpha\bigr) \hat{\rho}(\xi) + \mathrm{i} \xi
'\cdot\widehat{\beta\rho}(\xi)=0,\quad\quad \xi\in\RR^d,
\end{equation}
where $c_\alpha=\int_{\RR^d\backslash\{0\} }(\cos(u'y)-1)\nu
_\alpha(\dd y)$ for $u$ some unit vector in $\RR^d$ and $|\cdot|$
denotes the euclidean norm in $\RR^d$.
\end{corollary}

\begin{pf}
Let $\mu(\dd x)=\rho(x)\,\dd x$ be invariant for $X$. Since $\beta$ is
bounded and locally Lipschitz, and $Z$ and $W$ only act additively, we
can apply Proposition~\ref{proplevytwo} in the given setting and
obtain that
\begin{eqnarray*}
0&=& \int_{\RR^d} \mathrm{e}^{\mathrm{i}x'\xi} p(x,\xi) \rho(x)\,
\dd x
\\
&=& \int_{\RR^d} \mathrm{e}^{\mathrm{i}x'\xi} \bigl(-\mathrm{i}
\beta(x)'\xi+ a_1 |\xi|^2 - a_2
c_\alpha|\xi|^\alpha\bigr) \rho(x)\,\dd x
\\
&=& - \mathrm{i} \xi' \int_{\RR^d} \mathrm{e}^{\mathrm{i}x'\xi}
\beta(x) \rho(x)\,\dd x + \bigl(a_1 |\xi|^2 -
a_2 c_\alpha|\xi|^\alpha\bigr) \int
_{\RR^d} \mathrm {e}^{\mathrm{i}x\xi} \rho(x)\,\dd x.
\end{eqnarray*}
Substituting $\xi$ by $-\xi$ we further observe that this is
equivalent to
\[
0= \mathrm{i} \xi' \widehat{\beta\rho}(\xi) + \bigl(a_1 |
\xi|^2 - a_2 c_\alpha |\xi|^\alpha\bigr)
\hat{\rho}(\xi)
\]
which is \eqref{albecriterion}.
\end{pf}

\begin{remark}
In \cite{alberuediwu}, Proposition~3.1, an invariance condition for the type
of process considered in Corollary~\ref{cor:albeverio} is given.
Unfortunately, in their computations, the authors missed to use complex
conjugates when applying Parseval's identity, resulting in a wrong sign
(\cite{alberuediwu}, between equations (3.4) and (3.5)). The (corrected)
condition stated there is then automatically fulfilled if \eqref
{albecriterion} holds.
\end{remark}

Another interesting class of processes in our setting are processes
with factorizing symbol as they appear in the following corollary.

\begin{corollary}
Let $(X_t)_{t\geq0}$ be the unique solution of the SDE
\[
\dd X_t=\Phi(X_{t-})\,\dd L_t,\quad\quad t\geq0,
\]
where $\Phi\dvtx \RR\to\RR$ is bounded and locally Lipschitz continuous
and $(L_t)_{t\geq0}$ is a symmetric $\alpha$-stable, real-valued L\'
evy process, $\alpha\in(0,1)$.

Assume $X$ has an absolutely continuous invariant law $\mu(\dd x)=\rho
(x)\,\dd x$. Then $\Phi(x)\rho(x)= 0$ for $\lambda$-a.a. $x$.
\end{corollary}
\begin{pf}
We know from Proposition~\ref{proplevy} that the symbol of $X$ is
given by
$p(x, \xi)= |\Phi(x)|^\alpha|\xi|^\alpha$.

Suppose now that $X$ had an invariant law $\mu(\dd x)=\rho(x)\,\dd x$.
By our above results this yields
\[
0=\int\mathrm{e}^{\mathrm{i}x\xi} p(x,\xi) \mu(\dd x)= |\xi|^\alpha\int
\mathrm{e}^{\mathrm{i}x\xi} \bigl|\Phi(x)\bigr|^\alpha\rho(x) \,\dd x \quad\quad\forall\xi \in
\RR.
\]
Thus,
$f(\xi):=\int\mathrm{e}^{\mathrm{i}x\xi} |\Phi(x)|^\alpha\rho(x) \,\dd x=0$
for $\xi$ non-zero. Since $\Phi$ was assumed to be bounded, the
product $|\Phi(x)|^\alpha\rho(x)\leq C\cdot\rho(x)$ is integrable.
Hence, its Fourier transform is in $C_0$. Thus, $f(0)=0$ which gives
the claim.
\end{pf}

\subsection{Diffusions}

In case of a Brownian motion as driving process instead of a general L\'
evy process, we can even drop the boundedness condition on the
coefficient $\Phi$.

\begin{proposition} \label{corbrown}
Let $(X_t)_{t\geq0}$ be the unique solution of the SDE
\[
\dd X_t=-aX_t \,\dd t + \Phi(X_t) \,\dd
W_t,\quad\quad t\geq0,\vadjust{\goodbreak}
\]
where $(W_t)_{t\geq0}$ is an $n$-dimensional standard Brownian motion,
$a\in\RR$ and $\Phi\dvtx \RR^d\to\RR^{d\times n}$ is continuously
differentiable with bounded derivative.
Then $X$ is an It\^o process, for every $\xi\in\bbr^d$ the limit
$p(x,\xi)=\lim_{t\downarrow0} \lambda_\xi(x,t)$
exists and the functions $\lambda_\xi$ are globally bounded in $x$
(and $t$) for every $\xi\in\bbr^d$. Furthermore if $\mu$ is an
invariant law for $X$, then
\[
\int_{\RR^d} \mathrm{e}^{\mathrm{i}x' \xi} p(x,\xi) \mu(
\dd x)=\int_{\RR
^d} \mathrm{e}^{\mathrm{i}x' \xi} \bigl(\bigl|
\Phi(x)'\xi\bigr|^2+\mathrm{i}ax'\xi \bigr) \mu (\dd x)
=0 \quad\quad\mbox{for all } \xi\in\RR^d.
\]
\end{proposition}

We postpone the proof of this proposition to Section~\ref{Sproof}.

\begin{example} \label{exampleOU}
Consider the Ornstein--Uhlenbeck (OU) process driven by a
one-dimensional, standard Brownian motion $(W_t)_{t\geq0}$ with
parameters $\lambda>0$, $\sigma>0$
and starting random variable $X_0$, independent of $(W_t)_{t\geq0}$,
which is given by
\[
X_t=\mathrm{e}^{-\lambda t} \biggl( X_0 + \int
_{(0,t]} \mathrm {e}^{\lambda s} \sigma\,\dd W_s
\biggr) ,\quad\quad t\geq0.
\]
This process is a special case of the generalized OU process introduced
in Example~\ref{exampleGOUnec}. In particular, $X$ solves the SDE $\dd
X_t=-\lambda X_{t-}\,\dd t + \sigma\,\dd W_t$ such that we can now obtain
directly from Proposition~\ref{corbrown} that the symbol of the OU
process is
\[
p(x,\xi)=\mathrm{i}\lambda x\xi+ |\sigma\xi|^2.
\]
It is well known that $X$ admits a stationary distribution $\mu=\cL
(X_\infty)$ which is normal with mean
$0$ and variance $\frac{\sigma^2}{\lambda}$.
Using the symbol and this stationary distribution yields
\begin{eqnarray*}
\int_{\RR} \mathrm{e}^{\mathrm{i}x \xi} p(x,\xi) \mu(\dd x) &=&
\lambda\xi E\bigl[\mathrm{i}X_\infty\mathrm{e}^{\mathrm{i}X_\infty\xi}\bigr] + \sigma
^2\xi^2 E\bigl[\mathrm{e}^{\mathrm{i}X_\infty\xi}\bigr]
\\
&=& \lambda\xi\phi_{X_\infty}'(\xi)+ \sigma^2
\xi^2 \phi _{X_\infty}(\xi)
\\
&=& - \xi^2\sigma^2 \exp \biggl(-\frac{\sigma^2}{2\lambda} \xi
^2 \biggr) + \sigma^2\xi^2 \exp \biggl(-
\frac{\sigma^2}{2\lambda} \xi^2 \biggr)
\\
&=&0
\end{eqnarray*}
so that equation \eqref{symbolstat} is fulfilled.
\end{example}

\begin{example}
Let $(X_t)_{t\geq0}$ be the stochastic exponential of a Brownian motion
$(W_t)_{t\geq0}$ with variance $\sigma^2$, that is, $X_t=1+\int_{(0,t]}
X_{t-} \,\dd W_t$, then we have (\cite{protter}, Theorem II.37),
\[
X_t=\exp \bigl(W_t-\tfrac{1}{2}
\sigma^2 t \bigr),\quad\quad t\geq0.
\]
From Proposition~\ref{corbrown},
we obtain the corresponding symbol of the stochastic exponential as
\[
p(x,\xi)= x^2 \xi^2.
\]
Now, if $X$ had a stationary distribution $\mu=\cL(X_\infty)$ with
finite second moment, this would fulfill \eqref{symbolstat}, that is,
\[
0 = \int_{\RR} \mathrm{e}^{\mathrm{i}x \xi} x^2
\xi^2 \mu(\dd x)= \xi^2 \phi_{X_\infty}''(
\xi)
\]
and hence we had $\phi_{X_\infty}''(\xi)=0$ for all $\xi$ which is
only possible if $\phi_{X_\infty} =1$ for all $\xi$. Thus
$\mu$ had to be the Dirac measure at $0$. But obviously $X_t>0$ for
all $t\geq0$ which leads to a contradiction.
\end{example}

\section{Sufficiency}\label{S4}

As mentioned in the \hyperref[S1]{Introduction}, for Markov processes which are not
rich Feller, equation \eqref{generatorstat} is in general not
sufficient to prove invariance of the law $\mu$. Therefore we restrict
ourselves in this section to infinitesimal invariant laws, that is, to
laws which fulfill \eqref{generatorstat}.

\begin{theorem} \label{thmsuff} Let $(X_t)_{t\geq0}$ be an\vspace*{1.5pt} It\^o
process with generator $\cA$ whose domain $D(\cA)$ contains the test
functions $C_c^\infty(\bbr^d)$ and with symbol $p(x,\xi)$. Assume
there exists a probability measure $\mu$ such that $\int_{\RR^d} |
p(x,\xi)| \mu(\dd x)<\infty$ and $\int_{\RR^d} \mathrm{e}^{\mathrm{i}x'\xi
} p(x,\xi) \mu(\dd x)=0$. Then
\[
\int_{\RR^d} \cA f(x)\mu(\dd x)=0 \quad\quad\mbox{for all } f \in
C_c^\infty\bigl(\RR^d\bigr).
\]
\end{theorem}
\begin{pf}
By Lemma~\ref{itogenerator} the generator $\cA$ admits the
representation \eqref{generatorsymbol} for all $f\in C_c^\infty(\RR
^d)$. Hence using Fubini's theorem we obtain for all $f\in C_c^\infty
(\RR^d)$ 
\begin{eqnarray*}
0 &=& - \int_{\RR^d} \hat{f}(\xi) \int_{\RR^d}
\mathrm {e}^{\mathrm{i}x'\xi} p(x,\xi) \mu(\dd x) \,\dd\xi
\\
&=& - \int_{\RR^d} \int_{\RR^d}
\mathrm{e}^{\mathrm{i}x'\xi} p(x,\xi) \hat{f}(\xi) \,\dd\xi\mu(\dd x)
\\
&=&\int_{\RR^d} \cA f(x)\mu(\dd x).
\end{eqnarray*}
\upqed\end{pf}

The above theorem can easily be adapted to specific classes of symbols.
We illustrate this with the following corollary.

\begin{corollary}
Let $(X_t)_{t\geq0}$ be the unique solution of the SDE
\[
\dd X_t=-aX_t \,\dd t + \Phi(X_t) \,\dd
W_t,\quad\quad t\geq0,
\]
where $(W_t)_{t\geq0}$ is a standard Brownian motion, $a\in\RR$ and
$\Phi\dvtx \RR\to\RR$ is continuously differentiable with bounded
derivative and such that $|\Phi(x)|\leq K |x|^{\kappa/2}$ for some
constant $K$ and some $\kappa\in[1,2]$. Further suppose that the
domain $D(\cA)$ of the generator of $X$ contains the test functions
$C_c^\infty(\bbr^d)$. Assume there exists a probability distribution
$\mu$ such that $\int_{\RR^d} \mathrm{e}^{\mathrm{i}x'\xi} p(x,\xi) \mu
(\dd x)=0$ and $\int\|x\|^\kappa\mu(\dd x)<\infty$. Then
\[
\int_{\RR^d} \cA f(x)\mu(\dd x)=0 \quad\quad\mbox{for all } f \in
C_c^\infty\bigl(\RR^d\bigr).
\]
\end{corollary}
\begin{pf}
We know from Proposition~\ref{corbrown} that the symbol of $X$ is
given by
$p(x,\xi)= |\Phi(x)'\xi|^2+\mathrm{i}ax'\xi$. Hence, we have
$\int_{\RR^d} |\mathrm{e}^{\mathrm{i}x'\xi} p(x,\xi)| \mu(\dd x)\leq\int_{\RR^d} |p(x,\xi)| \mu(\dd x)<\infty$. By Theorem~\ref{thmsuff}
this gives the claim.
\end{pf}

\begin{example}
Let $(X_t)_{t\geq0}$ be a generalized Ornstein--Uhlenbeck process, as
defined in Example~\ref{exampleGOUnec}. Then by the same arguments as
in Example~\ref{exampleGOUnec} we see from Theorem~\ref{thmsuff}
together with \cite{liggett}, Theorem~3.37, that
\begin{equation}
\int_{\RR} \mathrm{e}^{\mathrm{i}x\xi} \psi_U(x
\xi) \mu (\dd x)=- \psi_L(\xi) \phi_\mu(\xi),\quad\quad \xi\in\RR,
\end{equation}
is also sufficient for $\mu$ to be an invariant law for $X$ with
finite second moment.

In the special case of the Ornstein--Uhlenbeck process as introduced in
Example~\ref{exampleOU} it is sufficient to suppose $\mu$ to be
integrable and equation \eqref{eqGOUstat} reduces to
\[
-\lambda\xi\phi_{\mu}'(\xi)= \sigma^2
\xi^2 \phi_{\mu}(\xi ),\quad\quad \xi\in\RR.
\]
This differential equation can be uniquely solved by $\phi_{\mu}(\xi
)=\exp(-\frac{\sigma^2}{2\lambda} \xi^2)$ (compare Example~\ref
{exampleOU}).
\end{example}

\section{Proofs} \label{Sproof}

\begin{pf*}{Proof of Lemma~\ref{glmkonv}}
We give the one-dimensional proof, since the multidimensional version
works alike; only the notation becomes more
involved. Let $x, \xi\in\bbr$. First, we use It\^o's formula under
the expectation and \def\theequation{\Roman{equation}}obtain
\begin{eqnarray}
\label{termone}\frac{1}{t} \bbe^x \bigl( \mathrm{e}^{\mathrm{i}(X_t-x)\xi} -1 \bigr)
&=& \frac{1}{t} \bbe^x \biggl(\int_{0+}^t
\mathrm{i} \xi\mathrm {e}^{\mathrm{i}(X_{s-}-x)\xi} \,\dd X_s \biggr)
\\
&&\label{termtwo}{}+ \frac{1}{t} \bbe^x \biggl(\frac{1}{2} \int
_{0+}^t-\xi^2 \mathrm
{e}^{\mathrm{i}(X_{s-}-x)\xi} \,\dd[X,X]_s^c \biggr)
\\
&&\label{termthree}{}+ \frac{1}{t} \bbe^x \biggl(\mathrm{e}^{-\mathrm{i}x\xi} \sum
_{0<s\leq t} \bigl(\mathrm{e}^{\mathrm{i} \xi X_s}-
\mathrm{e}^{\mathrm{i}\xi X_{s-}}-\mathrm{i}\xi\mathrm {e}^{\mathrm{i}\xi X_{s-}} \Delta X_s
\bigr) \biggr). 
\end{eqnarray}
In what follows, we will deal with the terms one-by-one.
To calculate term (\ref{termone}) we use the canonical decomposition of a
semimartingale (see \cite{jacodshir}, Theorem II.2.34) which we
write as follows\def\theequation{\arabic{section}.\arabic{equation}}\setcounter{equation}{0}
\begin{equation}
\label{candec} X_t=X_0 + X_t^{c}
+ \int_0^{t} \chi(y)y \bigl(\mu^{X}(
\cdot;\dd s,\dd y)-\nu(\cdot;\dd s,\dd y) \bigr) +\check{X}_t(
\chi)+B_t(\chi),
\end{equation}
where $\check{X}_t=\sum_{s\leq t}(\Delta X_s (1- \chi(\Delta X_s))$.
Therefore, term \eqref{termone} can be rewritten as
\begin{eqnarray*}
&&\frac{1}{t} \bbe^x \biggl(\int_{0+}^t
\mathrm{i} \xi\mathrm{e}^{\mathrm{i}(X_{s-}-x)
\xi} \,\dd \biggl( \underbrace{X_t^{c}}_{\mathrm{(IV)}}
 + \underbrace{\int_0^{t} \chi(y)y \bigl(
\mu^{X}(\cdot;\dd s,\dd y)-\nu(\cdot;\dd s,\dd y) \bigr)}_{\mathrm{(V)}}
\\
&&\hphantom{\frac{1}{t} \bbe^x \biggl(\int_{0+}^t
\mathrm{i} \xi\mathrm{e}^{\mathrm{i}(X_{s-}-x)
\xi} \,\dd \biggl(}{}+\underbrace {\check{X}_t(\chi)}_{\mathrm{(VI)}} +
\underbrace{B_t(\chi) }_{\mathrm{(VII)}} \biggr) \biggr).
\end{eqnarray*}
We use the linearity of the stochastic integral mapping. First, we
prove for term (IV)
\[
\bbe^x\int_{0+}^t \mathrm{i} \xi
\mathrm{e}^{\mathrm{i}(X_{s-}-x)\xi} \,\dd X_s^{c} = 0.
\]
The integral $\mathrm{e}^{\mathrm{i}(X_{t-}-x)\xi} \bullet X_t^{c}:=\int_{0+}^t \mathrm{e}^{\mathrm{i}(X_{s-}-x)\xi} \,\dd X_s^{c} $ is a local
martingale, since $X_t^{c}$ is a local martingale. To see that it is
indeed a martingale, we calculate
\[
\bigl[\mathrm{e}^{\mathrm{i}(X-x)\xi} \bullet X^{c}, \mathrm{e}^{\mathrm{i}(X-x)\xi}
\bullet X^{c} \bigr]_t = \int_0^t
\bigl(\mathrm{e}^{\mathrm{i}(X_s-x)\xi}\bigr)^2 \,\dd\bigl[X^{c},X^{c}
\bigr]_s = \int_0^t \bigl(\bigl(
\mathrm{e}^{\mathrm{i}(X_s-x)\xi}\bigr)^2 Q(X_s) \bigr) \,\dd s.
\]
The last term is uniformly bounded in $\omega$ and therefore, finite
for every $t\geq0$.
Hence, $\mathrm{e}^{\mathrm{i}(X_t-x)\xi} \bullet X_t^{c}$ is an
$L^2$-martingale which is zero at zero and therefore, its expected
value is constantly zero.

The same is true for the integrand (V): We show that the function
$H_{x,\xi}(\omega,s,y):=\mathrm{e}^{\mathrm{i}(X_{s-}-x)\xi}\cdot y \chi
(y)$ is in the class $F_p^2$ of Ikeda and Watanabe (see \cite{ikedawat}, Section
II), that is,
\[
\bbe^x \int_0^t \int
_{y\neq0} \bigl\llvert \mathrm{e}^{\mathrm{i}(X_{s-}-x)\xi} \cdot y\chi(y)
\bigr\rrvert ^2 \nu(\cdot;\dd s,\dd y) <\infty.
\]
In order to prove this, we observe
\[
\bbe^x \int_0^t \int
_{y\neq0} \bigl\llvert \mathrm{e}^{\mathrm{i}(X_{s-}-x)\xi
}\bigr\rrvert
^2 \cdot\bigl\llvert y\chi(y)\bigr\rrvert ^2 \nu(\cdot;
\dd s,\dd y) = \bbe^x \int_0^t \int
_{y\neq0} \bigl\llvert y\chi(y)\bigr\rrvert ^2
N(X_s,\dd y) \,\dd s.
\]
Since we have by hypothesis $\norm{\int_{y\neq0}(1\wedge y^2)
N(\cdot,\dd y)}_\infty< \infty$ this expected value is finite.
Therefore, the function $H_{x,\xi}$ is in $F_p^2$ and we conclude that
\begin{eqnarray*}
&&\int_0^t \mathrm{e}^{\mathrm{i}(X_{s-}-x)\xi} \,\dd
\biggl(\int_0^{s} \int_{y\neq0}
\chi(y)y \bigl(\mu^{X}(\cdot;dr,\dd y)-\nu(\cdot;dr,\dd y)\bigr)
\biggr)
\\
&&\quad=\int_0^t\int_{y\neq0}
\bigl(\mathrm {e}^{\mathrm{i}(X_{s-}-x)\xi} \chi(y)y \bigr) \bigl(\mu^{X}(\cdot;
\dd s,\dd y)-\nu(\cdot;\dd s,\dd y)\bigr)
\end{eqnarray*}
is a martingale. The last equality follows from \cite{jacodshir}, Theorem~I.1.30.

Now we deal with term \eqref{termtwo}. Here we have
\[
[X,X]_t^c= \bigl[X^c, X^c
\bigr]_t=C_t  =\bigl(Q(X_{t}) \bullet t\bigr)
\]
and therefore,
\begin{equation}
\label{liml} \frac{1}{2} \int_{0+}^t-
\xi^2 \mathrm{e}^{\mathrm{i}(X_{s-}-x)\xi} \,\dd[X,X]_s^c
 = - \frac{1}{2} \xi^2 \int_{0}^t
\mathrm{e}^{\mathrm{i}(X_{s-}-x)\xi} Q(X_{s}) \,\dd s.
\end{equation}
Since $Q$ is finely continuous and bounded we obtain by dominated convergence
\[
-\lim_{t\downarrow0} \frac{1}{2}\xi^2
\frac{1}{t} \bbe^x \int_0^t
\mathrm{e}^{\mathrm{i} (X_{s}-x)\xi} Q(X_{s}) \,\dd s =-\frac{1}{2}\xi
^2 Q(x).
\]
For the finite variation part of the first term, that is, (VII), we
obtain analogously
\begin{equation}
\label{limq} \lim_{t\downarrow0} \mathrm{i}\xi\frac{1}{t}
\bbe^x \int_0^t \mathrm
{e}^{\mathrm{i} (X_{s}-x)\xi} \ell(X_{s}) \,\dd s =\mathrm{i}\xi\ell(x).
\end{equation}
Finally, we have to deal with the various jump parts. At first, we
write the sum in \eqref{termthree} as an integral with respect to the
jump measure $\mu^{X}$ of the process:
\begin{eqnarray*}
&& \mathrm{e}^{-\mathrm{i}x\xi} \sum_{0<s\leq t} \bigl(
\mathrm{e}^{\mathrm{i}X_s\xi
}-\mathrm{e}^{\mathrm{i}X_{s-}\xi}-\mathrm{i}\xi\mathrm{e}^{\mathrm{i}\xi X_{s-}}
\Delta X_s \bigr)
\\
&&\quad= \mathrm{e}^{-\mathrm{i}x\xi} \sum_{0<s\leq t} \bigl(
\mathrm {e}^{\mathrm{i} X_{s-}\xi}\bigl(\mathrm{e}^{\mathrm{i}\xi\Delta X_s}-1 -\mathrm{i}\xi\Delta
X_s\bigr) \bigr)
\\
&&\quad=\int_{]0,t]\times\bbr^d} \bigl( \mathrm{e}^{\mathrm{i}
(X_{s-}-x)\xi}\bigl(
\mathrm{e}^{\mathrm{i}\xi y}-1 -\mathrm{i}\xi y\bigr)1_{\{y\neq0\}} \bigr)
\mu^{X}(\cdot;\dd s,\dd y)
\\
&&\quad=\int_{]0,t]\times\{y\neq0\}} \bigl( \mathrm {e}^{\mathrm{i} (X_{s-}-x) \xi}\bigl(
\mathrm{e}^{\mathrm{i}\xi y}-1 -\mathrm{i} \xi y\chi(y) -\mathrm{i}\xi y \cdot\bigl(1-\chi(y)\bigr)
\bigr) \bigr) \mu^{X}(\cdot;\dd s,\dd y)
\\
&&\quad= \int_{]0,t]\times\{y\neq0\}} \bigl(\mathrm {e}^{\mathrm{i} (X_{s-}-x) \xi}\bigl(
\mathrm{e}^{\mathrm{i}\xi y}-1 -\mathrm{i} \xi y\chi(y)\bigr) \bigr) \mu^{X}(\cdot;
\dd s,\dd y)
\\
&&\quad\quad{}+ \int_{]0,t]\times\{y\neq0\}} \bigl( \mathrm {e}^{\mathrm{i} (X_{s-}-x) \xi}\bigl( -\mathrm{i}\xi
y \cdot\bigl(1-\chi(y)\bigr)\bigr) \bigr) \mu ^{X}(\cdot;\dd s,\dd y).
\end{eqnarray*}
The last term cancels with the one we left behind from \eqref
{termone}, given by (VI). For the remainder-term, we get:
\begin{eqnarray*}
&& \frac{1}{t} \bbe^x \int_{]0,t]\times\{y\neq0\}} \bigl(
\mathrm{e}^{\mathrm{i} (X_{s-}-x) \xi}\bigl(\mathrm{e}^{\mathrm{i}\xi y}-1 -\mathrm{i} \xi y\chi(y)\bigr)
\bigr) \mu^{X}(\cdot;\dd s,\dd y)
\\
&&\quad= \frac{1}{t} \bbe^x \int_{]0,t]\times\{y\neq0\}} \bigl(
\mathrm{e}^{\mathrm{i} (X_{s-}-x) \xi}\bigl(\mathrm{e}^{\mathrm{i}\xi y}-1 -\mathrm{i} \xi y\chi(y)\bigr)
\bigr) \nu(\cdot;\dd s,\dd y)
\\
&&\quad= \frac{1}{t} \bbe^x \int_{]0,t]\times\{y\neq0\}} \bigl(
\underbrace{\mathrm{e}^{\mathrm{i} (X_{s-}-x) \xi}\bigl(\mathrm{e}^{\mathrm{i}\xi
y}-1 -\mathrm{i} \xi y
\chi(y)\bigr) \bigr)}_{:=g(s-,\cdot)} N(X_s,\dd y) \,\dd s.
\end{eqnarray*}

Here we have used the fact that it is possible to integrate with
respect to the compensator of a random measure instead of the measure
itself, if the integrand is in $F_p^1$ (see \cite{ikedawat}, Section
II.3). The function $g(s,\omega)$ is measurable and bounded
by our assumption, since $\abs{\mathrm{e}^{\mathrm{i}\xi y}-1 -\mathrm{i} \xi y\chi
(y)} \leq C_\xi\cdot(1\wedge\abs{y}^2)$, for a constant $C_\xi>0$.
Hence, $g\in F_p^1$.

Again by bounded convergence, we obtain
\begin{eqnarray}
\label{limn} &&\lim_{t \downarrow0} \frac{1}{t}
\bbe^x \int_{0}^t \mathrm
{e}^{\mathrm{i}(X_{s}-x) \xi} \int_{y\neq0} \bigl(\mathrm{e}^{\mathrm{i}y \xi}-1-\mathrm{i}y
\xi\chi(y) \bigr) N(X_{s},\dd y) \,\dd s\nonumber
\\[-8pt]\\[-8pt]
&&\quad= \int_{y\neq0} \bigl(\mathrm{e}^{\mathrm{i} y \xi}- 1 - \mathrm{i}y \xi
\chi(y) \bigr) N(x, \dd y).\nonumber
\end{eqnarray}
This is the last part of the symbol. Here we have used the continuity
assumption on $N(x,\dd y)$.

Considering the above calculations, in particular \eqref{liml}, \eqref
{limq} and \eqref{limn} we obtain
\begin{eqnarray*}
\biggl\llvert \frac{\bbe^x \mathrm{e}^{\mathrm{i}(X_t-x)'\xi}
-1}{t}\biggr\rrvert &=& \biggl|\mathrm{i}\xi\frac{1}{t}
\bbe^x \int_0^t
\mathrm{e}^{\mathrm{i} (X_{s-}-x)\xi
} \ell(X_{s}) \,\dd s - \frac{1}{2}
\xi^2\frac{1}{t} \bbe^x \int_{0}^t
\mathrm{e}^{\mathrm{i}(X_{s-}-x)\xi} Q(X_{s}) \,\dd s
\\
&&{}{}\,+ \frac{1}{t} \bbe^x \int_{0}^t
\mathrm{e}^{\mathrm{i}(X_{s-}-x) \xi} \int_{y\neq0} \bigl(
\mathrm{e}^{\mathrm{i}y \xi}-1-\mathrm{i}y \xi\chi(y) \bigr) N(X_{s},\dd y) \,\dd s
\biggr|
\\
&\leq&\llvert \xi\rrvert \frac{t}{t} \llVert \ell\rrVert _\infty+
\xi^2 \frac
{t}{2t} \llVert Q\rrVert _\infty+
C_\xi\frac{t}{t} \biggl\llVert \int_{y\neq0}
\bigl(1\wedge\llvert y\rrvert ^2\bigr) N(\cdot,\dd y) \biggr\rrVert
_\infty,
\end{eqnarray*}
a bound which is uniform in $t$ and $x$.
\end{pf*}

For the proof of Proposition~\ref{proplevy}, we need the following
lemma. Observe that for $\kappa\geq2$ and in the one-dimensional
case, this lemma follows directly from \cite{protter}, Theorem~V.67.

\begin{lemma} \label{lemma-levymoment}
Let $\kappa\geq1$ and suppose $(L_t)_{t\geq0}$ is a L\'evy process
such that $\EE[\|L_1\|^\kappa]<\infty$. Assume $X_0$ is a random
variable, independent of $L$, such that $\EE[\|X_0\|^\kappa]<\infty
$. Then the process $(X_t)_{t\geq0}$ defined by
\[
X_t=X_0 - a \int_{(0,t]}
X_{s-} \,\dd s + \int_{(0,t]} \Phi(X_{s-})
\,\dd L_s,
\]
where $\Phi\dvtx \RR^d\to\RR^{d\times n}$ is bounded, locally Lipschitz
and $a\in\RR$, fulfills
\[
\EE\Bigl[\sup_{0\leq t \leq1} \|X_t\|^\kappa\Bigr]<
\infty.
\]
\end{lemma}
\begin{pf}
Observe that
\[
\|X_t\|^\kappa 
\leq 4^\kappa
\|X_0\|^\kappa+ 4^\kappa|a|^\kappa\biggl
\llVert \int_{(0,t]} X_{s-} \,\dd s\biggr\rrVert
^\kappa+ 2^\kappa\biggl\llVert \int_{(0,t]}
\Phi(X_{s-}) \,\dd L_s\biggr\rrVert ^\kappa
\]
and hence for any $0\leq s\leq1$
\begin{eqnarray*}
&& \EE \Bigl[\sup_{0\leq t \leq s} \|X_t
\|^\kappa \Bigr]
\\
&&\quad\leq4^\kappa\EE\bigl[\|X_0\|^\kappa\bigr] +
4^\kappa|a|^\kappa\EE \biggl[\sup_{0\leq t \leq s} \biggl
\llVert \int_{(0,t]} X_{u} \,\dd u\biggr\rrVert
^\kappa \biggr] + 2^\kappa\EE \biggl[\sup_{0\leq t \leq s}
\biggl\llVert \int_{(0,t]} \Phi(X_{u-}) \,\dd
L_u\biggr\rrVert ^\kappa \biggr].
\end{eqnarray*}
By \cite{protter}, Lemma on bottom of page 345, we have that
\[
\EE \biggl[\sup_{0\leq t \leq s} \biggl\|\int_{(0,t]}
X_{u} \,\dd u\biggr\|^\kappa \biggr]\leq\int_{(0,s]}
E \bigl[\| X_{u}\|^\kappa \bigr] \,\dd u.
\]
On the other hand, it follows from an easy multivariate extension of
\cite{behme2011}, Lemma~6.1, that $\EE[\sup_{0\leq t \leq1} \|\int_{(0,t]} \Phi(X_{u-}) \,\dd L_u\|^\kappa]$ is finite, say $\leq K$,
under the given conditions. Thus
\[
\EE \Bigl[\sup_{0\leq t \leq s} \|X_t\|^\kappa
\Bigr] 
\leq 4^\kappa\EE \bigl[\|X_0
\|^\kappa \bigr] + 2^\kappa K + 4^\kappa|a|^\kappa
\int_{(0,s]} \EE \Bigl[\sup_{0\leq v \leq u}\|
X_{v}\|^\kappa \Bigr] \,\dd u.
\]
Now it follows from Gronwall's inequality (cf. \cite{protter}, Theorem
V.68) that
\[
\EE \Bigl[\sup_{0\leq t \leq1} \|X_t\|^\kappa
\Bigr]\leq\bigl(4^\kappa \EE \bigl[\|X_0\|^\kappa
\bigr] + 2^\kappa K\bigr) \mathrm{e}^{4^\kappa
|a|^\kappa}<\infty
\]
as we had to show.
\end{pf}

\begin{pf*}{Proof of Proposition~\ref{proplevy}}
It is well known that
the given SDE has a unique solution under the given conditions (cf.,
e.g., \cite{jacodshir}, Chapter~IX.6.7).
To keep notation simple, we give only the proof for $d=n=1$. Fix $x,\xi
\in\bbr$ and apply It\^o's formula to the function $\exp(\mathrm{i}(\cdot
-x)\xi)$:
\begin{eqnarray}
\label{threeterms} \frac{1}t \bbe^x \bigl( \mathrm{e}^{\mathrm{i}(X_t-x)\xi}
-1 \bigr)&=&\frac{1}t \bbe^x \biggl( \int
_{0+}^t \mathrm{i} \xi\mathrm{e}^{\mathrm{i} (X_{s-}-x) \xi} \,\dd
X_s - \frac{1}2 \int_{0+}^{t}
\xi^2 \mathrm{e}^{\mathrm{i} (X_{s-}-x) \xi} \,\dd[X,X]_s^c\nonumber
\\[-8pt]\\[-8pt]
&&\hphantom{\frac{1}t \bbe^x \biggl(}{}+\mathrm{e}^{-\mathrm{i}x \xi} \sum_{0<s \leq t} \bigl(
\mathrm {e}^{\mathrm{i} X_s \xi} - \mathrm{e}^{\mathrm{i} X_{s-} \xi} -\mathrm{i}\xi\mathrm
{e}^{\mathrm{i}X_{s-} \xi}\Delta X_s \bigr) \biggr).
\nonumber
\end{eqnarray}
For the first term, we get
\begin{eqnarray}
&&\frac{1}t \bbe^x  \int_{0+}^{t}
\bigl(\mathrm{i} \xi \mathrm{e}^{\mathrm{i}
(X_{s-}-x) \xi} \bigr) \,\dd X_s
\nonumber
\\
&&\quad= \frac{1}t \bbe^x \int_{0+}^{t}
\bigl(\mathrm{i} \xi \mathrm{e}^{\mathrm{i}
(X_{s-}-x) \xi} \bigr) \,\dd \biggl(\int
_0^s\Phi(X_{r-}) \,\dd
L_r \biggr)
\nonumber
- \frac{1}t \bbe^x \int
_{0+}^{t} \bigl(\mathrm{i} \xi \mathrm{e}^{\mathrm{i} (X_{s-}-x) \xi} a
X_{s-} \bigr) \,\dd s
\\
&&\quad= \frac{1}t \bbe^x \int_{0+}^{t}
\bigl(\mathrm{i} \xi \mathrm{e}^{\mathrm{i}
(X_{s-}-x) \xi} \Phi(X_{s-}) \bigr) \,\dd(\ell
s)\label{drift}
\\
&&\quad\quad{} + \frac{1}t \bbe^x \int_{0+}^{t}
\bigl(\mathrm{i} \xi \mathrm {e}^{\mathrm{i} (X_{s-}-x) \xi} \Phi(X_{s-}) \bigr) \,\dd
\biggl(\sum_{0<r\leq s} \Delta L_r
1_{\{\llvert \Delta Z_r\rrvert \geq1 \}} \biggr)\label{jumps}
\\
&&\quad\quad{} - \frac{1}t \bbe^x \int_{0+}^{t}
\bigl(\mathrm{i} \xi \mathrm {e}^{\mathrm{i} (X_{s-}-x) \xi} a X_{s-} \bigr) \,\dd s,
\label{driftzwei}
\end{eqnarray}
where we have used the L\'evy--It\^o decomposition of the L\'evy
process. Since the integrand is bounded, the martingale parts of the L\'
evy process yield martingales whose expected value is zero.

Now we deal with \eqref{jumps}. Adding this integral to the third
expression on the right-hand side of \eqref{threeterms} we obtain
\begin{eqnarray*}
&&\frac{1}t \bbe^x   \sum_{0<s \leq t}
\bigl(\mathrm {e}^{\mathrm{i}(X_{s-}-x)\xi} \bigl( \mathrm{e}^{\mathrm{i} \Phi(X_{s-})\Delta L_s \xi} -1- \mathrm{i}\xi\Phi
(X_{s-})\Delta L_s 1_{\{\llvert \Delta X_s\rrvert < 1
\}} \bigr) \bigr)
\\
& &\quad\accentset{t\downarrow0} {\longrightarrow} \int_{\bbr\setminus\{
0\}}
\bigl( \mathrm{e}^{\mathrm{i} \Phi(x)y \xi} -1- \mathrm{i}\xi\Phi(x)y 1_{\{\llvert y\rrvert < 1 \}} \bigr) N(\dd
y).
\end{eqnarray*}
The calculation above uses the same well-known results about
integration with respect to integer valued random measures as the proof
of Lemma~\ref{glmkonv}. In the case of a L\'evy process, the
compensator is of the form $\nu(\cdot;\dd s,\dd y)=N(\dd y) \,\dd s
$, see \cite{ikedawat}, Example II.4.2.

For the first drift part \eqref{drift}, we obtain
\[
\frac{1}t \bbe^x \int_{0+}^{t}
\bigl(\mathrm{i} \xi\cdot\mathrm{e}^{\mathrm{i}
(X_{s-}-x) \xi} \Phi(X_{s-})\ell \bigr) \,
\dd s = \mathrm{i} \xi\ell\cdot\bbe^x \frac{1}t\int
_{0}^{t} \bigl(\mathrm{e}^{\mathrm{i} (X_{s}-x) \xi}
\Phi(X_{s}) \bigr) \,\dd s \,
\accentset{t\downarrow0} {\longrightarrow}\, \mathrm{i}
\xi\ell\Phi(x).
\]

To deal with the second expression on the right-hand side of
\eqref{threeterms}, we first have to calculated the square bracket of
the process
\[
[X,X]_t^c = \biggl( \biggl[\int
_0^\cdot\Phi(X_{r-}) \,\dd
L_r, \int_0^\cdot
\Phi(X_{r-}) \,\dd L_r \biggr]_t^c
\biggr)=\int_0^t \Phi (X_{s-})^2
\,\dd(Q s).
\]
Let us remark that $\int a X_s \,\dd s$ is negligible in calculating
the square bracket $[X,X]_t$ since it is quadratic pure jump by \cite
{protter}, Theorem II.26. Now we can calculate the limit for the second
term of \eqref{threeterms}
\begin{eqnarray}\label{bracketlim}
&&\frac{1}{2t} \bbe^x \int_{0+}^{t}
\bigl(-\xi^2 \mathrm{e}^{\mathrm{i}
(X_{s-}-x) \xi} \bigr) \,\dd[X,X]_s^c\nonumber\\
&&\quad= \frac{1}{2t} \bbe^x \int_{0+}^{t}
\bigl(-\xi^2 \mathrm {e}^{\mathrm{i} (X_{s-}-x) \xi} \bigr) \,\dd \biggl(\int
_{0}^s \bigl(\Phi (X_{r-})
\bigr)^2 Q \,\dd r \biggr)
\nonumber
\\[-8pt]\\[-8pt]
&&\quad=-\frac{1}{2} \xi^2 Q \bbe^x \biggl(
\frac{1}t \int_{0}^t \bigl(
\mathrm{e}^{\mathrm{i}(X_s-x)\xi}\Phi(X_{s})^2 \,\dd s \bigr)
\biggr) \nonumber
\\
&&\quad\!\!\accentset{t\downarrow0} {\longrightarrow} -\frac{1}{2} \xi^2
Q \Phi(x)^2.
\nonumber
\end{eqnarray}
While in these three parts, due to the boundedness of $\Phi$, the
uniform boundedness of the approximants is trivially seen, we have to
be a bit more careful in dealing with the term \eqref{driftzwei}:
we use the Lemma~\ref{lemma-levymoment} and the fact $\sup_{0\leq t
\leq1} \bbe|X_t|\leq\bbe[\sup_{0\leq t \leq1} |X_t|]$. In order
to show that
\[
a\mathrm{i}\xi\bbe^x \int_0^1
\mathrm{e}^{\mathrm{i}(X_{ts}-x)\xi} X_{(ts)-} \,\dd s \underaccent{t\downarrow0} {
\longrightarrow} a\mathrm{i}\xi x
\]
in a uniformly bounded way, we consider
\begin{eqnarray*}
&&\bbe^x\int_0^1 \bigl\llvert
\mathrm{e}^{\mathrm{i}(X_{st}-x)\xi} X_{st}-\mathrm {e}^{\mathrm{i}(X_{st}-x)\xi} x +
\mathrm{e}^{\mathrm{i}(X_{st}-x)\xi} x -x\bigr\rrvert \,\dd s
\\
&&\quad= \bbe^x\int_0^1 \bigl\llvert
\mathrm{e}^{\mathrm{i}(X_{st}-x)\xi} (X_{st}-x) + \bigl(\mathrm{e}^{\mathrm{i}(X_{st}-x)\xi}-1
\bigr) x\bigr\rrvert \,\dd s.
\end{eqnarray*}
By $\bbe^x\abs{X_{st}-x}\leq c <\infty$, we can interchange the
order of integration.
In the end, we obtain
\begin{eqnarray*}
p(x,\xi) &=& -\mathrm{i} \ell\bigl(\Phi(x)\xi\bigr) +\mathrm{i}ax\xi+ \frac{1}{2} \bigl(
\Phi(x)\xi\bigr)Q \bigl(\Phi (x)\xi\bigr)
\\
&&{} -\int_{y\neq0} \bigl(\mathrm{e}^{\mathrm{i} (\Phi(x)\xi)y} -1 - \mathrm{i} \bigl(
\Phi(x)\xi\bigr)y \cdot1_{\{\llvert y\rrvert <1\}}(y) \bigr) N(\dd y)
\\
&=&\psi_L\bigl(\Phi(x)\xi\bigr)+\mathrm{i}ax\xi.
\end{eqnarray*}
Let us remark that in the multi-dimensional case the matrix $\Phi(x)$
has to be transposed, that is, the symbol of the solution is $\psi
_L(\Phi(x)'\xi)+\mathrm{i}ax'\xi$.

The result now follows as in the proof of Theorem~\ref{thm:nec}.
\end{pf*}

\begin{pf*}{Proof of Proposition~\ref{proplevytwo}}
In order to prove this result, we can mimic the previous proof. In this
case, $a=0$, the driving L\'evy process is $(Z',L')'\in\bbr^{d+n}$
and the bounded coefficient is $(b\cdot I_d,\Phi(x))\in\bbr^{d\times
(d+n)}$ where $I_d$ denotes the $d$-dimensional identity matrix. Since
$a$ is zero the respective part of the proof~-- the one where the moment
assumption is needed~-- can be omitted.
\end{pf*}

\begin{pf*}{Proof of Proposition~\ref{corbrown}}
The proof works perfectly analogue to the one of Proposition~\ref
{proplevy} with the following exception: from \cite{protter}, Theorem~V.67, we obtain that $\sup_{0\leq t \leq1} E(X_t)^2$ is
finite. This is needed in order to obtain the convergence in \eqref
{bracketlim} in a uniformly bounded way. In the present setting, $Q$ is
the identity matrix.
\end{pf*}


\section*{Acknowledgements}
Our thanks go to Ren\'e Schilling for helpful literature suggestions
and to the anonymous referee for his/her effort. Alexander Schnurr
gratefully acknowledges financial support by the German Science
Foundation (DFG) for the project SCHN1231/1-1 and the SFB 823 (project C5).


%

\printhistory
\end{document}